  \date{\today }   

 \documentclass[11pt,a4paper]{article}
 \usepackage{amsmath,amsfonts,amssymb,amsthm}
 \usepackage{enumerate,color}
\usepackage[applemac]{inputenc}
 \usepackage[polutonikogreek,english]{babel}

\newcommand{\ds}{\displaystyle}

\theoremstyle{plain}
\newtheorem{theorem}{Theorem}[section]
\newtheorem{proposition}[theorem]{Proposition}
\newtheorem{corollary}[theorem]{Corollary}
\newtheorem{lemma}[theorem]{Lemma}

\numberwithin{equation}{section}

\DeclareMathOperator{\conf}{conf}%

\DeclareMathOperator{\ccap}{Cap_{n}}
\DeclareMathOperator{\vol}{vol}%
\DeclareMathOperator{\lconf}{lconf}%

\renewcommand{\r}{\mathbb{R}}%
\newcommand{\tq}{\, \big| \, }%

\author{Vladimir Gol'dshtein and Marc Troyanov} 
\title{A Conformal de Rham Complex}
 
\begin{document}
 
\maketitle

\tableofcontents

\bigskip

\hrule

\begin{abstract}
We introduce the notion of a conformal de Rham complex of a Riemannian manifold. This is  a graded differential Banach algebra and it is invariant under quasiconformal  maps, in particular the associated cohomology is a new quasiconformal invariant. 

\medskip

\end{abstract}

\section{Introduction}

A conformal map between two domains or Riemannian manifolds is a diffeomorphism preserving all angles. A quasi-conformal map
can be thought of as a homeomorphism which sends infinitesimal spheres to infinitesimal ellipsoids w hoses axes are  uniformly bounded.
 The precise definition will be recalled below. The theory of quasi-conformal mappings between plane domains and Riemann surfaces goes back to the 1920s, it lies at the heart of Teuchmller theory, i.e. the theory of deformation of Riemann surfaces. It also plays an important role in elliptic partial differential equation and in some chapters of applied mathematics as illustrated for instance in the book \cite{bers} by Lipman Bers.

Lavrentiev suggested to study higher dimensional quasi-conformal mappings in the late 1930s, but systematic investigation began in the early 1960s by Yu.G. Reshetnyak and F. Gehring. In 1968, G. Mostow used the theory  quasi-conformal mappings in the proof of his celebrated rigidity theorem. Since then, the subject has continuously attracted the attention of many mathematicians. Recent developments include some generalization to subriemannian spaces such as Carnot groups, and more abstract metric measure spaces. See \cite{heinonen2006} for a very short description of the subject.

\medskip

One of the important question in the subject is to design \emph{quasi-conformal invariants}, i.e. invariants of Riemannian manifolds which
are stable under quasiconformal mappings and can thus be used to distinguish non quasiconformally equivalent manifolds.
An important class of quasiconformal invariants has been derived from the notion of \emph{conformal capacity} of condensers or, equivalently of \emph{moduli of curves families}. Another invariant is the so called \emph{Royden algebra}: if $M$ is a Riemannian manifold, we define its   
\emph{Royden algebra} $\mathcal{R}^n(M)$ to be the space of continuous functions $u : M \to \r$  such that $du \in L^n(M)$, where $n$ is the dimension of $M$. This is a Banach algebra for the norm
$$
 \|u \|_{\mathcal{R}^n(M)} = \|u \|_{L^{\infty}(M)} + \|u \|_{L^n(M)}.
$$
By the work of  M. Nakai, L. Lewis and J. Lelong-Ferrand, we know that the  Royden algebra is a complete quasiconformal invariant, that is $M$ and $N$ are quasiconformally equivalent manifolds if and only if $\mathcal{R}^n(M)$ and $\mathcal{R}^n(N)$ are isomorphic as abstract Banach algebra. Another example of algebraic quasiconformal invariant
is given by the  \emph{Dirichlet space} $\mathcal{L}^{1,n}(M)$ of all locally integrable functions $u : M \to \r$ such that 
$du \in L^n(M)$.  It is proved by S.K.Vodop'janov and V.Gol'dshtein in \cite{VG75} that two domains with Lipschitz boundaries  $U,V \subset \r^n$ are quasiconformally equivalent  if and only if there exists a lattice isomorphism between $\mathcal{L}^{1,n}(U)$ and $\mathcal{L}^{1,n}(V)$.
A third known algebraic invariant is the the space $BMO(M)$ of functions with bounded mean oscillation. It has been proved by Martin Reimann in \cite{reimann74} that two  quasiconformally equivalent manifolds have isomorphic $BMO$ spaces.
These algebraic invariants are important from a theoretical viewpoint, but one cannot use the Royden algebra or the Dirichlet space to quasiconformally distinguish two concrete manifolds, because these invariant are not  really computable.

\medskip

In the present paper, we describe a version of the de Rham complex adapted to quasiconformal  geometry. This is a Banach  differential graded algebra which is invariant under quasiconformal mappings by Theorem \ref{th.qcinvariance} below. We call this graded algebra the \emph{conformal de Rham complex} and denote it  by $\Omega_{\conf}^{\bullet}(M)$, it contains the Royden algebra in its center. We can then define an  associated cohomology: this \emph{conformal cohomology} is then obviously a quasiconformal invariant with  potentially interesting applications. Contrary to the Royden algebra, it is not completely hopeless to  try and compute this conformal cohomology; we give in the paper some partial results in this direction. As an application we prove that the 3-dimensional Lie group SOL is not quasiconformaly equivalent to the hyperbolic space $\mathbb{H}^{3}$. This fact is based on a computation of the
cohomology in degree 2 and could not be proved by standard methods based on capacities, because both spaces are homogeneous spaces with exponential volume growth.

\medskip 

The paper is organized as follow. We first recall some basic facts on the notion of weak exterior derivative for locally integrable differential forms. We then define the conformal de Rham complex $\Omega_{\conf}^{\bullet}(M)$ in section 3 and we prove that it is a Banach differential algebra. In section 4, we prove a chain rule for forms in $\Omega_{\conf}^{\bullet}$ and maps of class $W^{1,n}$. Section 5 is a digression on conformal capacities and the quasiconformal invariance of $\Omega_{\conf}^{\bullet}$ is proved in section 6. Section 7  contains definitions and results on the cohomology of  the conformal de Rham complex. Section 8 studies a notion of
\emph{interior}  conformal de Rham complex, containing a notion of differential forms ``vanishing at infinity'' and section 9 studies the top dimensional cohomology. The last section contains some applications and some remarks on related subjects.

\section{Locally integrable differential forms}

Recall that if  $(M,g)$ is a Riemannian manifold and $x\in M$, then the $k^{th}$ exterior power of the cotangent space $\Lambda^k T_x^*M$ inherits a scalar product defined by
$$
 \langle \theta^1 \wedge \theta^2 \wedge \cdots \wedge \theta^k , 
 \varphi^1 \wedge  \varphi^2 \wedge \cdots \wedge  \varphi^k \rangle = \det (g(\theta^i , \varphi^j)),
$$ 
where \ $\theta^i , \varphi^j \in \Lambda^1 T_x^*M$. 

A \emph{measurable differential form} $\theta$ of degree $k$ on $M$ is a measurable section of the vector bundle
$\Lambda^k T_x^*M \to M$. Such a form $\theta$ is said to be  \emph{locally integrable} if 
$$
  \int_U \| \theta \| d\vol_g < \infty
$$
for every relatively compact subset $U \subset\subset M$. We denote by $L_{loc}^1(M,\Lambda^k)$ the space of locally integrable differential forms of degree $k$ on $M$, and by $C^r(M,\Lambda^k)$
the space of  $k$-forms of class $C^r$. Finally $C^r_0(M,\Lambda^k)$ is the space of $k$-forms of class $C^r$ with compact support.

\medskip

The exterior differential is a well known operator $d : C^r(M,\Lambda^k) \to C^{r-1}(M,\Lambda^{k+1})$. There is also  a notion of weak exterior differential for elements in $L_{loc}^1(M,\Lambda^k)$:

\medskip

\textbf{Definitions} 
\textbf{a.)} One says that a form  $\theta \in L^1_{loc}(M,\Lambda^{k+1})$ is the
\emph{weak exterior differential}  (or the exterior differential \emph{in the sense of currents}) of a form
$\phi \in L^1_{loc}(M,\Lambda^{k})$ and one writes
{$d\phi = \theta$} if one has
\[
 \int_U \theta \wedge \omega = (-1)^{k+1}\int_U \phi \wedge d\omega
 \,  .
\]
for any oriented open subset $U\subset M$ and any smooth form 
$\omega \in C_0^{\infty}(U,\Lambda^{n-k})$
 with compact support in $U$.

\medskip

\textbf{b.)}  A sequence $\{\alpha_j \} \subset L^1_{loc}(M, \Lambda^k)$ is said to 
\emph{converge weakly}  to $\alpha \in  L^1_{loc}(M, \Lambda^k)$ if and only if 
$$\lim_{j \to \infty}\int_U \alpha_j\wedge \omega   \to \int_U \alpha\wedge \omega \, .$$
for any oriented open subset $U\subset M$ and any $\omega \in C_0^{\infty}(U,\Lambda^{n-k})$.
Convergence in $L^1_{loc}$ implies weak convergence.

\medskip

\begin{proposition} \label{prop.1}
 The weak exterior differential satisfies the following properties:
 \begin{enumerate}[i)]
  \item If $\theta \in L^1_{loc}(M,\Lambda^k)$ satisfies $\theta =d\phi $ for some $\phi \in L^1_{loc}(M,\Lambda^{k-1})$, then $d\theta = 0$.
  \item Let $\alpha \in L^1_{loc}(M, \Lambda^k)$ and
$\beta \in  L^1_{loc}(M, \Lambda^{k+1})$.
If there exists a sequence $\{\alpha_j \} \subset  C^1(M, \Lambda^k)$ such
that
$\alpha_j \to \alpha$ and $d\alpha_j \to \beta$ weakly, 
 then $d\alpha = \beta$.
\end{enumerate}
\end{proposition}

\medskip

\textbf{Proof}  (i) The first assertion is clear,  because for any smooth test form $\omega$, we have  $d d\omega = 0$. \\
(ii) For any  smooth $(n-k-1)-$form $\omega$  with compact support in an oriented domain $U \subset M$, we have from Stokes formula 
$
 \int_U d(\alpha\wedge \omega ) = 0.
$
Hence 
$$\int_U \alpha\wedge d\omega = \lim_{j\to \infty} \int_U \alpha_j\wedge d\omega
=  \lim_{j\to \infty} (-1)^{k+1}\int_U d\alpha_j\wedge \omega = 
(-1)^{k+1}\int_U \beta\wedge \omega \, .$$
By definition  of the weak exterior differential, this means that $d\alpha = \beta$.

\qed

\section{The conformal de Rham complex}

Given a Riemannian manifold $(M,g)$, we introduce the space 
$$
  \Omega_{\conf}^{k}(M,g) = \{ \omega \in L^{n/k}(M,\Lambda^k) \tq d\omega \in L^{n/(k+1)}(M,\Lambda^{k+1}) \}
$$
of differential forms of degree $k$ in  $L^{n/k}(M)$ having a weak exterior differential in $L^{n/(k+1)}$. It is a 
Banach space for the graph norm
\begin{equation} \label{confnorm}
 \| \omega \|_{\conf} =   \left( \int_M |\omega |^{n/k}d\vol_g\right)^{\frac{k}{n}}
+ \left( \int_M |d\omega|^{n/(k+1)}d\vol_g\right)^{\frac{k+1}{n}} .
\end{equation}

\bigskip

\textbf{Definition} The \emph{conformal de Rham complex} of a Riemannian manifold $(M,g)$ is the pair
$(\Omega_{\conf}^{\bullet}(M,g) , d)$, where $\Omega_{\conf}^{\bullet}(M,g) =\oplus_{k=1}^n \Omega_{\conf}^{k}(M,g)$.

\medskip

The name is justified by the following lemma:
\begin{lemma}\label{lem.invconf}
$(\Omega_{\conf}^{\bullet}(M,g) , d)$ is a conformal invariant
of $(M,g)$, that is if $g' = \lambda^2 \, g$ is a conformal deformation of the metric $g$ (where $\lambda$ is a smooth positive function on $M$), then $(\Omega_{\conf}^{\bullet}(M,g') , d) = (\Omega_{\conf}^{\bullet}(M,g) , d)$, and the norm  $\| \cdot \|_{\conf} $ is the same for both metrics.
\end{lemma}

\textbf{Proof}  A direct calculation shows that the $L^{n/k}$-norm is conformally invariant on $k$-forms, see e.g. 
\cite{GoT2006}.

\qed

\bigskip

 \textbf{Remark}   Observe  that $\Omega_{\conf}^{n}(M)$ is simply the space of $L^{1}(M,\Lambda^n)$ of integrable $n$-forms on $M$ and that  $\Omega_{\conf}^{0}(M)$ is the space of functions $u\in L^{\infty}(M)$ such that $du\in L^n(M)$.  The spaces $L^1$ and $L^{\infty}$ are not very well behaved Banach spaces, this causes some intricacies in our study of the  conformal de Rham complex.

\begin{lemma}\label{lem.density}
Smooth forms are dense in $\Omega_{\conf}^{k}(M)$ for $k>0$, that is 
  $$
   \overline{\Omega_{\conf}^{k}(M) \cap C^{\infty}(M, \Lambda^{k})} = 
   \Omega_{\conf}^{k}(M).
  $$
  Smooth functions are dense  in $\Omega_{\conf}^{0}(M)$  for the weak convergence defined above.
\end{lemma}

\textbf{Proof} This is proved by regularization, see e.g.  \cite{Burenkov,GT}, see also \cite{GoT2006}.

\qed

\bigskip

\textbf{Definition} A \emph{Banach differential graded algebra} is a triple $(A,\cdot , d)$ where 
\begin{enumerate}[a)]
  \item $A$ is a Banach space;
  \item $A$ is a direct sum $A=\oplus_{k\in \mathbb{N}}A^k$ where $A^k\subset A$ is a closed subspace;
  \item $d : A\to A$ is a bounded operator such that $d(A^k) \subset A^{k+1}$ and $d\circ d = 0$;
  \item $(A,\cdot )$ is a real (or complex) algebra such that 
  $x \in A^k, y \in A^l \Rightarrow x\cdot y \in A^{k+l}$;
  \item $A$ is commutative in the graded sense, that is 
  $x\cdot y =  (-1)^{kl}y\cdot  x $ for any $y\in A^l$ and $x\in A^k$;
  \item $(A,\cdot )$ is a Banach algebra, that is $\|x\cdot y\| \leq \|x\| \| y\|$;
  \item The Leibniz rule holds: $$d(x\cdot y) = d(x) \cdot y + (-1)^k x \cdot d(y),$$ for any $y\in A$ and $x\in A^k$ .
\end{enumerate}

\bigskip

\begin{theorem}
 The conformal de Rham complex $(\Omega_{\conf}^{\bullet}(M),  d)$ is a unitary Banach differential graded algebra
 for the exterior product $\wedge$.
\end{theorem}

\textbf{Proof} It is clear that    each \  $\Omega_{\conf}^{k}(M)$ is a Banach space and 
$d : \Omega_{\conf}^{k}(M) \to \Omega_{\conf}^{k+1}(M)$ is a bounded operator.
That $d\circ d = 0$ follows from Proposition \ref{prop.1}(i) above. 

We need to prove that  $\Omega_{\conf}^{\bullet}(M)$ is a {Banach algebra}, that is
 if $\alpha,\beta \in \Omega_{\conf}^{\bullet}(M) $, then 
 $\alpha\wedge\beta \in \Omega_{\conf}^{\bullet}(M)$ and
\begin{equation}\label{eq.normprod}
 \| \alpha\wedge\beta\|_{\conf} \leq \| \alpha\|_{\conf} \cdot \| \beta\|_{\conf}.
\end{equation}
Observe first  that
if $\alpha \in \Omega_{\conf}^{k}(M)$ and $\beta \in \Omega_{\conf}^{\ell}(M)$, then 
\begin{equation}\label{leibniz}
   \| \alpha\wedge\beta\|_{L^{n/(k+\ell)}} \leq \| \alpha \|_{L^{n/k}} \|\beta\|_{L^{n/\ell}}
\end{equation}
by Hlder's inequality, since $\frac{k}{n}+\frac{\ell}{n}=\frac{k+\ell}{n}$.  \ 
Likewise
\begin{equation*}
  \| (d\alpha)\wedge\beta\|_{L^{n/(k+\ell+1)}}  \leq  \| d\alpha \|_{L^{n/(k+1)}} \|\beta\|_{L^{n/\ell}}
 \end{equation*}
and
\begin{equation*}
 \|\alpha\wedge d(\beta)\|_{L^{n/(k+\ell+1)}}   \leq   \| \alpha\|_{L^{n/k}} \|d\beta\|_{L^{n/(\ell+1)}}.
 \end{equation*}
We   now prove (\ref{eq.normprod}) assuming $\alpha, \beta$ smooth. For such form the 
Leibniz rule  is known:
\begin{equation}
   d(\alpha\wedge\beta) = d(\alpha)\wedge\beta + (-1)^k\alpha\wedge (d\beta),
\end{equation}
and we have, using the  three inequalities above,
\begin{eqnarray*}
   \| \alpha\wedge\beta\|_{\conf}   & = &   
   \| \alpha\wedge\beta\|_{L^{n/(k+\ell)}} +   \| d(\alpha\wedge\beta)\|_{L^{n/(k+\ell +1)}} 
 \\ & = & 
    \| \alpha\wedge\beta\|_{L^{n/(k+\ell)}} +   \| d(\alpha)\wedge\beta \ \pm \ \alpha\wedge d(\beta)\|_{L^{n/(k+\ell +1)}} 
 \\ & \leq & 
     \| \alpha\wedge\beta\|_{L^{n/(k+\ell)}} +   \| d(\alpha)\wedge\beta\|_{L^{n/(k+\ell +1)}} +
     \| \alpha\wedge d(\beta)\|_{L^{n/(k+\ell +1)}} 
 \\ & \leq & 
  \| \alpha \|_{L^{n/k}} \|\beta\|_{L^{n/\ell}} +  \| d\alpha \|_{L^{n/(k+1)}} \|\beta\|_{L^{n/\ell}}
  +    \| \alpha\|_{L^{n/k}} \|d\beta\|_{L^{n/(\ell+1)}}
\\ & \leq & 
\left(  \| \alpha \|_{L^{n/k}}  +  \| d\alpha \|_{L^{n/(k+1)}}  \right) \cdot 
\left(  \|\beta\|_{L^{n/\ell}} +  \| d\beta\|_{L^{n/(\ell+1)}} \right)
\\ & = & 
 \| \alpha\|_{\conf} \cdot \| \beta\|_{\conf}.
\end{eqnarray*}
 By density of smooth forms for $k, \ell>0$, the inequality 
  (\ref{eq.normprod}) follows now for any forms in $\Omega_{\conf}^{\bullet}(M)$ as well as the Leibniz rule
  (\ref{leibniz}). If $k=0$ or $\ell = 0$, but in this case we only have sequences of forms
  weakly converging to $\alpha$ and $\beta$, the argument also work in this case by the lower semi-continuity of the 
  norm  $ \| \cdot  \|_{\conf}$ for the weak convergence.
  \\
  We have proved that  $\Omega_{\conf}^{\bullet}(M)$ is a  Banach differential graded algebra. It is 
  unitary, since $1\in \Omega_{\conf}^{0}(M) \subset L^{\infty}(M)$.

  \qed
  
  \bigskip
  
Recall  that  $\Omega_{\conf}^{0}(M)$  is a subalgebra of  $\Omega_{\conf}^{\bullet}(M)$. It is a closed subspace and thus
a Banach algebra with norm
$$\|u\|_{\conf} = \|u\|_{L^{\infty}(M)} + \|du\|_{L^n(M)}.$$
  
 \bigskip
   
\textbf{Definition} The \emph{conformal Royden algebra} of $M$ is the set
$$\mathcal{R}^n(M) = C(M) \cap \Omega_{\conf}^{0}(M)$$
of continuous function in $\Omega_{\conf}^{0}(M)$. The closure  of continuous functions with compact support in $\Omega_{\conf}^{0}(M)$ is  denoted  by  $\mathcal{R}_0^n(M)$. These are  closed subalgebras and we have the following sequence of closed subalgebras:
$$
 \mathcal{R}_0^n(M) \subset \mathcal{R}^n(M) \subset \Omega_{\conf}^{0}(M) \subset 
 \Omega_{\conf}^{\bullet}(M).
$$

 \bigskip

The Royden algebra is \emph{not} dense in 
 $\Omega_{\conf}^{0}(M)$, in particular, Lemma \ref{lem.density}  does not hold for $k=0$. 
 
 \medskip

\section{The chain rule for $W^{1,n}$-Sobolev maps.}

The following result shows the usefulness of the conformal de Rham complex.

\begin{proposition}\label{prop.chainrule1}
  Let $U \subset \r^n$, $V\subset \r^m$ be bounded domains and $f\in W^{1,n}(U,  V)$. Then for any smooth  
  differential form $\beta$   defined  on a neighborhood of  $V$, we have  \ $f^*\beta \in  \Omega_{\conf}^{\bullet}(U)$
  and
  $$
   f^*(d\beta) = d (f^*\beta).
  $$
\end{proposition}

\bigskip

\textbf{Proof}  \  As in the case of smooth mappings, this result follows from the chain rule for $0$-forms (i.e. functions) and the Leibniz rule. 

More precisely, let us write $f(x) = (f^1(x),f^2(x),\cdots f^m(x))$, by hypothesis, we have $f^j$ bounded and  $df^j \in L^n(U,\Lambda^1)$ for any $1 \leq j \leq m$, we thus  have $f^j \in 
 \Omega_{\conf}^{0}(U)$. 
 In particular, we have $df^j \in 
 \Omega_{\conf}^{1}(U)$ and since $ \Omega_{\conf}^{\bullet}(U)$ is an algebra, we have
\begin{equation}\label{regdf1}
 df^{i_1}\wedge df^{i_2}\wedge \cdots \wedge df^{i_k} \in  \Omega_{\conf}^{k}(U)
\end{equation}
for any   $1 \leq i_1,i_2,...,i_k \leq n$  (this is of course simply the H\"older inequality).

\medskip   

We first prove the Proposition for $0$-forms. So let us consider   an arbitrary smooth function $b$  defined  on a neighborhood of  $V$. The chain rule for functions and Sobolev maps (see e.g. \cite{GT}) implies  that 
$f^*(db) = d(b\circ f)$ in the weak sense. 
Observe that
$$
 d(b\circ f) = f^*(db) =  \sum_{\nu=1}^n \left(\frac{\partial b}{\partial y^{\nu}} \circ f \right)  df^{\nu} .
$$
Our hypothesis imply that  $b\circ f$ and $\frac{\partial b}{\partial y^{\nu}} \circ f$ are  bounded, hence we have $ d(b\circ f) = f^*(db)   \in L^n(U,\Lambda^1)$ and thus $b\circ f \in  \Omega_{\conf}^{0}(U)$.

\medskip

We can now prove the Proposition for a general $k$-form. By linearity, it is sufficient to 
consider  differential forms $\beta \in C^{\infty} (\r^n, \Lambda^k)$ of type 
$\beta = b(y) \,  dy^{i_1}\wedge dy^{i_2}\wedge \cdots \wedge dy^{i_k}$. 
We just proved that 
 $b\circ f$ and $df^j$ belong to $f^*\beta \in  \Omega_{\conf}^{\bullet}(U)$, thus 
$$
  f^*(\beta) =  (b\circ f) \  df^{i_1}\wedge df^{i_2}\wedge \cdots \wedge df^{i_k} 
   \in  \Omega_{\conf}^{\bullet}(U)
$$
because $ \Omega_{\conf}^{\bullet}(U)$ is an algebra. Now the  Leibniz rule, together with $d(b\circ f) = f^*(db)$ and $d(df^j) = 0$, implies that
\begin{eqnarray*}
  d(f^*(\beta) )& = & d\left( (b\circ f) \  df^{i_1}\wedge df^{i_2}\wedge \cdots \wedge df^{i_k} \right)
    \\ & = & 
     d (b\circ f) \wedge  df^{i_1}\wedge df^{i_2}\wedge \cdots \wedge df^{i_k}
      \\ & = & 
      f^*(db)\wedge  df^{i_1}\wedge df^{i_2}\wedge \cdots \wedge df^{i_k}
            \\ & = &   f^*(d\beta).
\end{eqnarray*}

\qed

\bigskip

\section{Conformal capacity}

\textbf{Definition} Let $M$ be $n$-dimensional Riemannian manifold and $F \subset U \subset M$ a pair of subsets with $U$ open in $M$.
The \emph{conformal capacity} of the pair \ $(F,U)$ is then defined as follow
\begin{equation*}
\ccap(F,U) := \inf \left\{ \int_U |du|^nd\vol \tq  u \in  \mathcal{A}(F,U)  \,  \right\} \, ,
\end{equation*}
where the set of admissible functions is given by
$$
  \mathcal{A}(F,U):= \{   u \in \mathcal{R}_0^n(M)  \tq   u\ge 1 \text{ on a neighbourhood of $F$ and } u \geq 0 \ \text{a.e.} \}
$$
If $\mathcal{A}_{p}(F,U) = \emptyset$, then we set
$\ccap(F,U)=\infty$.
If $U=M$, we simply write $\ccap (F,M)=\ccap(F)$.

\bigskip

There is a number of important notions related to the concept of capacity, in particular the notion of parabolic manifolds and that of polar sets.

\bigskip

\textbf{Definition}  A set \ $S\subset M$ \ is  \emph{conformally polar} if for any pair of open  relatively compact sets $U_{1} \subset U_{2} \subset M$ such that $\mathrm{dist}(U_{1},M\setminus U_2 )>0$, we have \  $\ccap(S\cap U_{1},U_{2})=0$.

\bigskip

Any finite set in $M$ is conformally polar while no set of Hausdorff dimension $>0$ is conformally polar. In particular
a conformally polar set is always totally discontinuous, it even has Hausdorff dimension $0$. See e.g. \cite{resh89} for these facts and more on polar sets.

\begin{theorem}\label{th.removable}
 Conformally polar sets are removable sets for the conformal de Rham complex, that is
$$\Omega_{\conf}^{\bullet}(M\setminus S) = \Omega_{\conf}^{\bullet}(M)$$
for any conformally polar subset $S \setminus M$.
\end{theorem} 

\bigskip

For instance $\Omega_{\conf}^{\bullet}(\r^n) = \Omega_{\conf}^{\bullet}(S^n)$, since the  Euclidean space is conformally equivalent to a sphere with a point removed.

\bigskip

\textbf{Proof}   Using Proposition 3.1 in \cite{tr99}, one can find a sequence
of function $u_j \in \Omega_{\conf}^{0}(M)$ such that $u_j = 0$ is a neighbourhhod of $S$ and $u_j \to 1$ uniformly
on any compact subset of $M \setminus S$ and $\int_M |du_j|^n \to 0$.
Using Lebesgue dominated convergence theorem, we see that
$$\lim_{j\to 0} \|(u_j\omega) - \omega \|_{L^{n/k}} =0 
\quad \mathrm{and }\quad 
\lim_{j\to 0} \|(u_jd\omega) - d\omega \|_{L^{n/(k+1)}} =0$$
for any $k$-form $\omega \in \Omega_{\conf}^{\bullet}(M)$. We also have from H\"older's inequality
$$
 \lim_{j\to 0} \|du_j \wedge \omega\|_{L^{n/(k+1)}}  \leq 
  \lim_{j\to \infty}  \|du_j\|_{L^{n}} \cdot  
  \|\omega\|_{L^{n/k}} \, = 0,
$$
thus, by the Leibniz rule, we have
$$
 \lim_{j\to 0} \|(u_j\omega) - \omega \|_{\conf} =0.
$$
It follows that the sequence of bounded operators
$$
 T_j : \Omega_{\conf}^{\bullet}(M\setminus S) \to \Omega_{\conf}^{\bullet}(M)
$$
defined by $T_j (\omega) = u_j \omega$ converges to an isometry between these
Banach algebras.

\qed

\bigskip

\textbf{Definitions}  The Riemannian manifold $M$ is \emph{conformally parabolic} if $\ccap(F,M)= 0$ for any compact set $F \subset M$. It is  \emph{conformally hyperbolic} otherwise. 

\bigskip

\begin{lemma}\label{lem.critparab1}
An $n$-dimensional Riemannian manifold $M$ is  conformal  parabolic if and only if there exists a sequence
 of  smooth functions with compact support $\{\eta_j\} \subset \mathcal{R}_0^n(M)$ such that $\eta_j \to 1$ uniformly on each compact set and   $\ds \lim_{j\to \infty}\int_M |d\eta_j|^n =0$. 
\end{lemma}

\bigskip  

\textbf{Proof}  Assume that  $M$ is conformally parabolic and choose a sequence $D_1 \subset D_2 \subset \cdots \subset M$ of compact sets whose union is $M$. By hypothesis, we have $\ccap(D_j, M) = 0$,
hence one can find for each $j$ a function $u_j \in \mathcal{R}_0^n(M)$  such that  $u_j \geq  1$ on $D_j$
and $\int_M |du_j|^n \leq 1/j$. The desired sequence is obtained by trunction: $\eta_j = \min \{ \max \{u_j, 1\} \}$.

\medskip

Conversely, if such a sequence  $\{\eta_j\} $ exists, then $\ccap(F,M) = 0$ for all compact set $F$ in $M$ by definition  and $M$ is thus conformally parabolic.
 
\qed

\bigskip

\textbf{Remark } The  notions of conformal capacity, parabolic manifolds and polar sets have been intensively studied, see e.g. \cite{GR,holopainen90,HKM,resh89,tr99,ZK} among other works. The conformal capacity is essentially equivalent to the notion of modulus of families of paths.

\section{Quasiconformal maps}

Quasiconformal maps  have been intensively studied since, see e.g. 
 \cite{vaisala79,vuorinen92,vuorinen2007,heinonen2006} for some background on this subject.
Let us recall that a map $f : (M,g) \to (N,h)$ between two $n$-dimensional Riemannian manifolds is said
to be \emph{quasi-conformal}, if it satisfies the following properties:
\begin{enumerate}[i.)]
  \item $f$ is a homeomorphism;
  \item $f \in W^{1,n}_{loc}(M,N)$, where $n= \dim (M) = \dim (N)$;
  \item there exists a constant $K$ such that 
\begin{equation}\label{ineqK}
  |df_x|^n \leq K \, |J_f(x)|
\end{equation}
almost everywhere, where $|df_x|$ is the norm of the weak differential of $f$ at $x\in M$ and $J_f(x)$ is its Jacobian.
\end{enumerate}

\medskip

This is the analytic definition, there is also a well known geometric definition:
\begin{theorem}
The  homeomorphism $f : (M,g) \to (N,h)$  is  quasiconformal if  and only if for all $x$ in $M$  we have
$$ \limsup_{r\to 0} \frac{\sup\{d(f(x),f(y)); d(x,y)=r\}} {\inf\{d(f(x),f(y)); d(x,y)=r\}}\le H $$ 
for some $H < \infty$. 
\end{theorem}
See e.g. \cite{vaisala79,vuorinen92}.

\qed

\bigskip

\textbf{Remarks}  A continuous map $f : M \to N$ satisfying the conditions (ii) and (iii) is called a \emph{quasiregular map} or a \emph{mapping with bounded distortion}. These maps are an important generalization of the class of quasiconformal maps, basic references on  them are the books \cite{resh89,rickman93}.

\medskip

\begin{theorem} Any quasiconformal $f : (M,g) \to (N,h)$ satisfies the following properties:
\begin{enumerate}[a)]
  \item $f$ is differentiable almost everywhere;
  \item $f$ maps sets of measure zero in $M$ to sets of mesure zero in $N$ (this is called the Luzin property);
  \item the inverse map $f^{-1} : N \to M$ is also quasi-conformal;
  \item the composition of two quasiconformal  maps is again a quasiconformal  map;
  \item the ``change of variables formula'' for integrals holds: for any  measurable function
  $v : N \to \r_+$ on $N$, the pull-back  $v\circ f$ is measurable on $M$ and
$$\int_M  (v\circ f) (x)|J_f(x) | \, d\vol_g(x) = 
\int_N    v(y)  \, d\vol_h(y)\, .$$
\end{enumerate}
\end{theorem}

\medskip

see  \cite{vaisala79,resh89} for  proofs of these facts.

\bigskip

We quote two additional classical results on quasiconformal maps:

\begin{theorem}
 Two manifolds $M$ and $N$ are quasiconformally equivalent if and only if the Royden algebras $\mathcal{R}^n_0(M)$  and $\mathcal{R}^n_0(N)$ are isomorphic as abstract Banach algebras.
\end{theorem}

\medskip

This is a deep result of  J. Lelong-Ferrand, see \cite{ferrand}.

\begin{theorem}\label{th.qcmapcapa}
 The homeomorphism $f : (M,g) \to (N,h)$ is quasiconformal if and only if for any pair of sets 
 $F\subset U \subset M$  with $U$ open and $F$ compact we have
 $$
   \frac{1}{K}\ccap{(F,U)} \leq   \ccap{(f(F),f(U))} \leq K \ccap{(F,U)}. 
 $$
\end{theorem}

\medskip

(See \cite{vaisala79} for a proof using   modulus of path families instead of capacities.)

\medskip

\begin{corollary}
 Let  $f : (M,g) \to (N,h)$ be a quasiconformal homeomorphism, then 
 \begin{enumerate}[i)]
  \item a subset $S\subset M$ is conformaly polar if and only if $f(S)\subset N$ is conformaly polar;
  \item $M$ is conformally parabolic if and only if $N$ is conformally parabolic.
\end{enumerate}
\end{corollary}

The proof is obvious after the previous Theorem.

\bigskip

\begin{theorem}\label{th.qcinvariance}
Let $f : (M,g) \to (N,h)$  be a homeomorphism between two Riemannian manifolds. Then $f$ is  a
quasiconformal map if and only if the pull-back of differential forms defines an isomorphism of Banach 
differential algebras
 $$
  f^*  : \Omega_{\conf}^{\bullet}(N) {\longrightarrow} \Omega_{\conf}^{\bullet}(M).
 $$
\end{theorem}

\bigskip

\textbf{Proof}  Assume that $f$ is a $K$-quasiconformal map. Then for any differential  form $\theta$ on $N$ of degree $k$, we have a.e.
$$|f^* \theta_x | \leq |df(x)|^k \, |\theta_{f(x)} | \leq K^{k/n} 
\,  |J_f(x)|^{k/n} \, |\theta_{f(x)} |
$$
where $K$ is the constant in inequality (\ref{ineqK}).
By the   change of variable formula, we then have 
$$
\int_M |f^*\theta_{f(x)} |^{n/k} \, d\vol_g  \leq K 
 \, \int_M |\theta_{f(x)} |^{n/k}  \, |J_f(x)| \,  d\vol_g \leq
K  \, \int_N |\theta_{f(x)} |^q \,   d\vol_h \, .
$$
This shows that
$\ds  \|f^*\theta\|_{L^{n/k}(M,\Lambda^k)} \leq K^{k/n}  \|\theta\|_{L^{n/k}(N,\Lambda^k)}$, likewise, 
we also have 
$$\ds  \|f^*(d\theta)\|_{L^{n/(k+1)}(M,\Lambda^{k+1})} \leq
 K^{(k+1)/n}  \|d\theta\|_{L^{n/(k+1)}(N,\Lambda^{k+1})}.$$
  We thus have proved that $f^*$ acts as a bounded operator on the de Rham complexes and 
$$
 \|f^*\theta\|_{\conf} \leq K^{(k+1)/n}  \|\theta\|_{\conf}.
$$
It remain to prove the chain rule $d(f^*\theta) = f^*(d\theta)$ for any form $\theta$ in the conformal de Rham complex. But for smooth forms, the chain rule follows from Proposition \ref{prop.chainrule1}. Since smooth forms are dense in the  conformal  de Rham complex and $ f^*  : \Omega_{\conf}^{\bullet}(N) \to \Omega_{\conf}^{\bullet}(M)$ is continuous, the chain rule holds for any $\theta \in  \Omega_{\conf}^{\bullet}(N)$.

\medskip

Conversely,  $f^* : C_0(N) \to C_0(M)$ is an isomorphism because $f$ is a homeomorphism.
If $f^*  : \Omega_{\conf}^{\bullet}(N) \overset{\simeq}{\longrightarrow} \Omega_{\conf}^{\bullet}(M)$ is also an isomorphism, then it defines an isomrphism at the level of the Royden algebras: $ f^*  : \mathcal{R}^n_0(N) \overset{\simeq}{\longrightarrow} \mathcal{R}^n_0(M)$. By definition, the conformal capacities in $M$ and $N$ are then comparable
and thus $f$ is quasiconformal by Theorem \ref{th.qcmapcapa}.

  \qed

\section{Conformal cohomology}

\subsection{Definitions}

Given an arbitrary Banach complex, one can define a \emph{reduced cohomology}, a \emph{non reduced cohomology} and a \emph{torsion}  (see \cite{GoT2006,Grom,pansu06}).  Let us recall these definitions in the present situation. First  set
$Z^k_{\conf}(M) :=L^{n/k}(M,\Lambda^k) \cap \ker d$
(it is the set of closed forms in $L^{n/k}(M,\Lambda^k)$) 
and
\[
  B^k_{\conf}(M) := d\left(L^{n/{k-1}}(M,\Lambda^{k-1}) \right) \cap
  L^{n/k}(M,\Lambda^k).
\]
Observe that $Z^k_{\conf}(M) \subset \Omega_{\conf}^{k}(U)$  is a closed linear subspace and that 
 $$ B^k_{\conf}(M) \subset \overline{B}^k_{\conf}(M) \subset Z^k_{\conf}(M). $$

\bigskip

\textbf{Definition}
The \emph{conformal cohomology} of $M$ is the quotient
 \begin{equation*}
H_{\conf}^{k}(M):=Z_{\conf}^{k}(M)/B_{\conf}^{k}(M)\,,
\end{equation*}
 the \emph{ reduced conformal cohomology} of $M$ is
 \begin{equation*}
\overline{H}_{\conf}^{k}(M):=Z_{\conf}^{k}(M)/\overline{B}_{\conf}^{k}(M)\,,
\end{equation*}
            
 (where $\overline{B}_{\conf}^{k}(M)$ is the closure of
$B_{\conf}^{k}(M)$) and the \emph{torsion} is
$$
  T_{\conf}^{k}(M):=\overline{B}_{\conf}^{k}(M)\,/ B_{\conf}^{k}(M).
$$
 
\textbf{Remarks:}  These are vector spaces. The reduced cohomology inherits a norm from
that of $Z_{\conf}^{k}(M)$ and it is a Banach space.  The unreduced
cohomology is a Banach space if and only if the torsion vanishes. If the  torsion $T_{\conf}^{k}(M)\neq 0$,
then it is infinite dimensional and we always have  the exact sequence
\begin{equation}\label{exact}
0\rightarrow T_{\conf}^{k}(M)\rightarrow H_{\conf}^{k}(M)\rightarrow
\overline{H}_{\conf}^{k}(M)\rightarrow 0.
\end{equation}
The conformal cohomology depends on the metric $g$, and if necessary, we denote it by 
$H_{\conf}^{k}(M,g)$ (and likewise for $\overline{H}_{\conf}^{k}(M,g)$ and $T_{\conf}^{k}(M,g)$.
However, it follows from lemma \ref{lem.invconf}  that these spaces are conformal invariant, they only depend on the conformal class of the metric.

\bigskip

\begin{lemma}                      
$Z^{\bullet}_{\conf}(M)$ is a subalgebra of  $\Omega^{\bullet}_{\conf}(M)$, while
$B^{\bullet}_{\conf}(M)$ and $\overline{B}^{\bullet}_{\conf}(M)$ are ideals in 
$Z^{\bullet}_{\conf}(M)$.                                                                                               
\end{lemma}

\bigskip  

\textbf{Proof}  The Leibniz rule  implies that if $\alpha,\beta \in \Omega^{\bullet}_{\conf}(M)$
are closed, then $d(\alpha \wedge \beta) = 0$, thus $Z^{\bullet}_{\conf}(M)$ is a subalgebra. Now if 
$\alpha = d\gamma \in B^{\bullet}_{\conf}(M)$ and $\beta \in Z^{\bullet}_{\conf}(M)$, then 
$$d(\gamma \wedge \beta) = \alpha \wedge \beta , $$
hence $\alpha \wedge \beta  \in B^{\bullet}_{\conf}(M)$ and that set is thus an ideal in $Z^{\bullet}_{\conf}(M)$. Now the closure of an ideal in a Banach algebra is clearly also an ideal. This concludes the proof.

\qed

\bigskip

This lemma immediately implies the following
\begin{corollary}
 The wedge product defines structures of algebra on $H_{\conf}^{\bullet}(\mathbb{B}^{n})$,
 $\overline{H}_{\conf}^{\bullet}(\mathbb{B}^{n})$ and $T_{\conf}^{\bullet}(\mathbb{B}^{n})$.
\end{corollary}

\bigskip

Observe also  the  exact sequence (\ref{exact}) holds in the category of graded algebras.

\subsection{Relation with Sobolev inequalities}

The
$L_{q,p}$-cohomology of a Riemannian manifold  $(M,g)$ has an interpretation
in terms of Sobolev inequalities for 
differential forms \cite{GoT2006}. In the conformal case, the results can be stated as follow:
\begin{theorem} \label{th.sobin1a}
\begin{enumerate}[(i)]
  \item $H_{\conf}^{k}(M,g)=0$
 if and only if there exists a constant $C<\infty$ such that for
 any closed form \\ 
  \ $\omega \in L^{n/k}(M,\Lambda^k)$
 there exists a differential form $\theta \in L^{n/(k-1)}(M,\Lambda^{k-1})$,
 such that $d\theta=\omega$ and
\[
\left\Vert \theta\right\Vert _{L^{k/(k-1)}} \leq C\left\Vert \omega\right\Vert _{L^{k/n}}.
\]
  \item If $T_{\conf}^{k}(M)=0$, \,
 then there exists a constant $C$ such that for any differential
form  $\theta \in L^{n/(n-k)}(M,\Lambda^{k-1})$ there exists a
form $\zeta\in L^{n/(n-k)}(M,\Lambda^{k-1})$ such that $d\zeta = 0$ and
\begin{equation}\label{inconc.sob2}
\left\Vert \theta-\zeta\right\Vert _{L^{k/(k-1)}}  \leq C\left\Vert
d\theta\right\Vert _{L^{n/k}}.
\end{equation}
  \item  Conversely, if $1<k<n$ and  if a Sobolev inequality 
(\ref{inconc.sob2}) holds, then $T_{\conf}^{k}(M)=0$.
\end{enumerate}
\end{theorem}

\bigskip
 
\begin{corollary}\label{cor.T1not0}
 For every manifolds, we have $T_{\conf}^{1}(M) \neq 0$.
\end{corollary}

\bigskip

\textbf{Proof}  Let $x\in M$ be an arbitrary point  and $U\subset M$ be a neighborhood of  $x$. Then we can find
an essentially unbounded function $u\in L^1(M)$ and such that $du \in L^n(M)$ (see  \cite[page 119]{stein}) . In particular, inequality (\ref{inconc.sob2})
never holds.

\qed 

\bigskip

\subsection{The Poincar\'e Lemma}

We have the following Poincar\'e Lemma:
\begin{theorem}[The Poincar\'e Lemma]\label{th.Poincarelemma}
$H_{\conf}^{k}(\mathbb{B}^{n})=0$
if   $1 < k  < n$. We also have 
$H_{\conf}^{0}(\mathbb{B}^{n}) = \overline{H}_{\conf}^{0}(\mathbb{B}^{n})=\mathbb{R}$ and
$H_{\conf}^{1}(\mathbb{B}^{n}) \neq 0$.
\end{theorem}

\textbf{Proof}  This is theorem 11.5 in  \cite{GoT2006} (in that paper, one should assume $q<\infty$ and not $q\leq \infty$). 
The fact that $H_{\conf}^{1}(\mathbb{B}^{n}) \neq 0$ follows from
Corollary \ref{cor.T1not0}, and the equality $H_{\conf}^{0}(\mathbb{B}^{n}) = \r$ is obvious.

\qed

\medskip

Whether  $H_{\conf}^{n}(\mathbb{B}^{n})$ and $\overline{H}_{\conf}^{n}(\mathbb{B}^{n})$ vanish is still an open question for us.

\bigskip

\begin{corollary}\label{cor.cohhyp}
 The conformal cohomology of the hyperbolic space $\mathbb{H}^n$ vanishes for any $1<k<n$.
\end{corollary}

\textbf{Proof}  The hyperbolic space is conformally equivalent to the unit ball via the Poincar ball model. 
Thus $H_{\conf}^{k}(\mathbb{H}^{n})=H_{\conf}^{k}(\mathbb{B}^{n})$.

\qed

\subsection{Other Results}

For conformally parabolic $n$-manifolds, we have a better result:
 
\begin{theorem}\label{computeHn1}
 Let $M$ be a connected oriented conformally parabolic $n$-manifolds
 then \\ 
 $\overline{H}_{\conf}^{n}(M)  \simeq \r$.
\end{theorem}

\medskip

This result is proved in Theorem \ref{computeHn2}  below.

\begin{theorem}
 Any cohomology class (or reduced cohomology class) of degree $1 < k < n$ in $M$ can be represented
 by a smooth form. 
\end{theorem}

\medskip

This is proved by regularisation, see Theorem 12.8 in \cite{GoT2006}. 

\qed

\bigskip

\begin{theorem}\label{th.qcinvariance2}
Any quasiconformal map $f : (M,g) \to (N,h)$  defines an isomorphism at the level of the
conformal cohomolgy  
 $$
  f^*  : H_{\conf}^{\bullet}(N) {\longrightarrow} H_{\conf}^{\bullet}(M),
 $$
 and likewise for reduced cohomology and torsion.
\end{theorem}

\textbf{Proof} This follows immediatey from  Theorem \ref{th.qcinvariance}.

\qed

\bigskip

\begin{theorem}\label{th.removability2}
Let $M$ be an arbitrary Riemannian manifold and $S\subset M$ a conformally polar set. Then 
$M\setminus S$ and $M$ have the same conformal cohomology and torsion.
\end{theorem}

\textbf{Proof}  It is immediate from  Theorem \ref{th.removable}.

\qed
 
\bigskip

\begin{corollary}
If $M$is compact and $N$ is quasiconformally equivalent to $M\setminus S$ where
$S\subset M$ is a conformally polar set, then 
$$
   H_{\conf}^{k}(N)=H_{de Rham}^{k}(M)
$$
and $T_{\conf}^{k}(N)=0$ for any $1 < k <n$.
\end{corollary} 

\bigskip

In particular, since $\r^n$ is conformally a sphere with a point removed, we have
$H_{\conf}^{k}(\r^n) = H^{k}(S^n) = 0$ for any $1 < k <n$.
We also have $H_{\conf}^{n}(\r^n) = H^{n}(S^n) = \r$  by Theorem  \ref{computeHn1}
and $H_{\conf}^{0}(\r^n) = H^{0}(S^n) = \r$ for obvious reasons. On the other hand,
$T_{\conf}^{1}(\r^n) \neq 0$ by  Corollary \ref{cor.T1not0}.

\section{Interior  conformal de Rham complex}

Let us denote by \ ${\Omega}_{\conf,0}^{k}(M)$ the closure of
smooth
forms with compact support in $\Omega _{\conf}^{k}(M)$ for the graph norm (\ref{confnorm}). 
Concretely, a $k$-form $\theta$ belongs to ${\Omega}_{\conf,0}^{k}(M)$ if and only if
$\theta \in {\Omega}_{\conf}^{k}(M)$ and there exists a sequence $\{\theta_j\}$ of smooth $k$-forms with compact support
 such that $d\theta_j \in L^{n/(k+1)}$ for all $j$ and 
$$
 \lim_{j\to \infty} \|\theta - \theta_j \|_{\conf} = 0.
$$

We also define $Z_{\conf,0}^{k}(M) = \ker(d) \cap {\Omega}_{\conf,0}^{k}(M)$ and 
$B_{\conf,0}^{k}(M) =  d({\Omega}_{\conf,0}^{k-1}(M))$. Observe that 
the closure $\overline{B}_{\conf,0}^{k}(M)$ is simply the closure 
of $d\left(C^{\infty}_0(M,\Lambda^k) \right)$
in the space $L^{n/k}(M,\Lambda^k)$.
The corresponding cohomologies are defined in the usual way:
$$
 H_{\conf,0}^{k}(M) = Z_{\conf,0}^{k}(M) / B_{\conf,0}^{k}(M) \quad \mathrm{and} \quad
  \overline{H}_{\conf,0}^{k}(M) = Z_{\conf,0}^{k}(M) / \overline{B}_{\conf,0}^{k}(M).
$$

\bigskip

\begin{proposition}\label{prop.interior} 
The interior conformal de Rham complex satisfies the following properties:
\begin{enumerate}[i.)]
  \item   ${\Omega}_{\conf,0}^{\bullet}(M) \subset {\Omega}_{\conf}^{\bullet}(M)$ is a closed ideal;
  \item  ${\Omega}_{\conf,0}^{0}(M) = \mathcal{R}^n_0(M)$.
\end{enumerate}
\end{proposition}

The proof is elementary.
  
\bigskip

\begin{theorem}\label{th.conf0isconf} 
If  $M$ is conformally parabolic then 
 ${\Omega}_{\conf,0}^{k}(M) = {\Omega}_{\conf}^{k}(M)$
 for any $k \geq 1$.
\end{theorem}

\bigskip

\textbf{Proof}   The proof is similar to that of Theorem  \ref{th.removable}: Choose a sequence $\{\eta_j\} \in \mathcal{R}^n_0(M)$ as in Lemma \ref{lem.critparab1}. It is easy to check that for any  $ \omega\in  \Omega_{\conf}^{k}(M)$, we have $\eta_j\omega \in  \Omega_{\conf,0}^{k}(M)$, and 
Lebesgue dominated convergence theorem implies that
$$
\lim_{j\to 0} \|(\eta_j\omega) - \omega \|_{L^{n/k}} =0 
\quad \mathrm{and }\quad 
\lim_{j\to 0} \|(\eta_jd\omega) - d\omega \|_{L^{n/(k+1)}} =0
$$
We also have from H\"older's inequality
$$
 \lim_{j\to 0} \|d\eta_j \wedge \omega\|_{L^{n/(k+1)}}  \leq 
  \lim_{j\to \infty}  \|d\eta_j\|_{L^{n}} \cdot  
  \|\omega\|_{L^{n/k}} \, = 0,
$$
thus, by the Leibniz rule, we have
$$
 \lim_{j\to 0} \|(\eta_j\omega) - \omega \|_{\conf} =0.
$$
It follows that  $\ds \omega = \lim_{j\to \infty} (\eta_j\omega) \in  \Omega_{\conf,0}^{k}(M)$,
and thus  ${\Omega}_{\conf,0}^{k}(M) = {\Omega}_{\conf}^{k}(M)$.

\qed

\bigskip

\textbf{Remark} The proof clearly does not work for $k=0$. In fact   we have 
$1 \in {\Omega}_{\conf}^{0}(M) \setminus \Omega_{\conf,0}^{0}(M)$ 
 for any non compact manifold (because for any $u \in C^{\infty}_0(M)$, we have
$\|1-u \|_{\conf} \geq \|1-u \|_{L^{\infty}} \geq 1$).

\section{The integration operator}

In this section, $M$ is an oriented $n$-dimensional Riemannian manifold.

Since   $\Omega_{\conf}^{n}(M)= L^{1}(M,\Lambda^n)$, we have a well defined  bounded linear form
\begin{equation}\label{int.op}
   I : \Omega_{\conf}^{n}(M) \to \r
\end{equation}
given by integration: $I(\omega) = \int_M\omega$.

 \bigskip

\begin{proposition} \label{prop.intoncohom}
Let $M$ be an oriented  Riemannian manifold of dimension $n$. Then the integration operator given by
 (\ref{int.op}) vanishes on  $\overline{B}_{\conf,0}^{n}(M)$.
 Hence there is is a well defined bounded linear form
 $$
 \overline{H}_{\conf,0}^{n}(M)  \to  \r 
 $$
 given by  $[\omega] \mapsto \int_M\omega$. This linear form is surjective, in particular
 $\overline{H}_{\conf,0}^{n}(M) \neq 0$.
\end{proposition}

\bigskip
  
\textbf{Proof} By Stokes theorem, the integration operator 
 vanishes on smooth exact forms with compact support and thus on  ${B}_{\conf,0}^{n}(M)$ and   $\overline{B}_{\conf,0}^{n}(M)$
 by density and continuity. It is thus well defined on the quotient 
 $\overline{H}_{\conf}^{n}(M) = \Omega_{\conf}^{n}(M) /\overline{B}_{\conf}^{n}(M)$.
Since we can always find a smooth $n$-form with compact support and non vanishing integral, this linear form is surjective.

\qed

\bigskip

For conformally parabolic $n$-manifolds, we have a better result:
 
\begin{theorem}\label{computeHn2}
 Let $M$ be a connected oriented conformally parabolic $n$-manifolds
 then the integration operator
 \begin{equation}\label{int.op2}
   I : \overline{H}_{\conf}^{n}(M)  \overset{\simeq}{\longrightarrow} \r
\end{equation}
defines an isomorphism.
\end{theorem}

\bigskip  

\textbf{Proof}   Since $M$ is conformally parabolic, we know from Theorem
\ref{th.conf0isconf} that  $\overline{H}_{\conf}^{n}(M)
= \overline{H}_{\conf, 0}^{n}(M)$ and thus the linear form $I$ is well defined. As observed above,
this linear form is surjective.

\medskip

By Theorem 3 in \cite{GoT99}, we know that a smooth $n$ form $\omega$  with compact support
has a primitive $\theta \in L^{n/(n-1)}$ if and only if $I(\omega) = 0$. In other words, if
$M$ is conformally parabolic, then 
$$
 \ker(I) \cap C^{\infty}_0(M,\Lambda^n) \subset B^n_{\conf}(M).
$$
But $C^{\infty}_0(M,\Lambda^n)$ is dense in $\Omega^n_{\conf,0}(M)=\Omega^n_{\conf}(M)$,  we thus
have the exact sequence
$$
 0 \to \overline{B}^n_{\conf}(M) \to \Omega^n_{\conf}(M) \overset{I}{\longrightarrow} \r \to 0.
$$
Therefore the operator  $I : \overline{H}_{\conf}^{n}(M) = \Omega_{\conf}^{n}(M) /\overline{B}_{\conf}^{n}(M) \to \r$ is an isomorphism.

\qed

\bigskip

In  contrast to the previous theorem, we have
  
\begin{theorem}[The Kelvin-Nevanlinna-Royden criterion] \label{th.KRN}
Let $M$ be an  oriented  Riemannian manifold of dimension $n$. Then 
$$
 I :  B_{\conf}^{n}(M) \to  \r
$$
is  non trivial if and only if $M$ is conformally hyperbolic. 
\end{theorem}

\medskip

This result implies in particular that an integration operator
$  I : \overline{H}_{\conf}^{n}(M)  \to \r$ 
is  well defined only  if $M$ is conformally parabolic.

\medskip

\textbf{Proof} If $M$ is conformally parabolic, then $I(B_{\conf}^{n}(M)) = I(B_{\conf,0}^{n}(M)) = 0$.
The converse direction is proved in \cite[Proposition 1]{GoT99},  which says that if $M$ is 
conformally hyperbolic, then there exists $\beta \in \Omega_{\conf}^{n-1}(M)$ such that $\int_M d\beta \neq 0$.

\qed
 
\section{Complements}

 
\subsection{An Application to SOL geometry}

Let us denote by  $\mathbb{H}^{3}$ the three dimensional hyperbolic space and by 
SOL  the group of $3\times 3$ real matrices of the form
$$ {\footnotesize
\left(
\begin{array}{ccc}
e^z & 0 & x \\
0 & e^{-z} & y \\
0 & 0 & 1
\end{array}
\right) \, . }$$
This  is a solvable and unimodular three dimensional Lie group. It is diffeomorphic to
$\r^{3}$ (with coordinates $x,y,z$) and a left invariant Riemannian
metric is $ds^2 = e^{-2z}dx^2 + e^{2z}dy^2 + dz^2$ its volume
measure is given by $dxdydz$ and is bi-invariant.  See \cite{tr98,SOL} for more information on the geometry of
this group.

\bigskip

\begin{theorem}
 $\mathbb{H}^{3}$ and SOL are not quasiconformally equivalent.
\end{theorem}

\bigskip

This is result is not easy to prove directly, because  $SOL$ and $\mathbb{H}^{3}$ are both
conformally hyperbolic, and they both have exponential volume growth.

\bigskip

\textbf{Proof}  $\mathbb{H}^{3}$ is conformally equivalent to the three ball  $\mathbb{B}^{3}$, thus
$H_{\conf}^{2}(\mathbb{H}^{3}) = H_{\conf}^{2}(\mathbb{B}^{3}) =0$. On the other hand, we have 
$H_{\conf}^{2}(SOL) \neq 0$, a fact proven in \cite{SOL}. The result follows then from Theorem 
\ref{th.qcinvariance}.

\qed

\subsection{Relative conformal cohomology}

Let us now set
\begin{equation*}
\hat{\Omega}_{\conf}^{k}(M) = \Omega _{\conf}^{k}(M)/{\Omega}_{\conf,0}^{k}(M),
\end{equation*}
the cohomology of this complex is the  \emph{relative conformal cohomology}, or  the 
\emph{conformal cohomology at infinity} of $M $, denoted by ${\hat{H}}_{\conf,0}^{k}(M)$. From the
short exact sequence of complexes
\begin{equation*}
0\rightarrow {\Omega}_{\conf,0}^{k}(M)\rightarrow \Omega _{\conf
}^{k}(M)\rightarrow \hat{\Omega} _{\conf}^{k}(M)\rightarrow 0,
\end{equation*}
we deduce a long exact sequence in cohomology
\begin{equation*} 
\cdots \rightarrow {H}_{\conf,0}^{k}(M)\rightarrow H _{\conf
}^{k}(M)\rightarrow \hat{H} _{\conf}^{k}(M)\rightarrow
{H}_{\conf,0}^{k+1}(M)\rightarrow \cdots .
\end{equation*}

\bigskip

Of course, this is useful only for conformally hyperbolic manifolds, since  we know from Theorem
\ref{th.conf0isconf} that  for a  conformally parabolic manifold
we have  $\hat{\Omega}_{k}^{\bullet}(M) = 0$ for all $k \geq 1$.

\subsection{The local conformal de Rham complex}

It also useful to consider a \emph{local  conformal de Rham complex}: we denote it by 
$\Omega_{\lconf}^{n}(M)$. By definition, a locally integrable differential form $\alpha$
belongs to $\Omega_{\lconf}^{k}(M)$ if and only if $h\cdot \alpha \in \Omega_{\conf}^{k}(M)$
for any smooth function $h : M \to \r$ with compact support. This local complex is a differential algebra,
but not a Banach space. It has been 
introduced in the work of  Sullivan and Donaldson, see \cite{DS}.  For any $k\neq 1,n$, the 
cohomology of this complex coincides with de Rham cohomology  (use the Poincar lemma 
(Theorem \ref{th.Poincarelemma}) and the  sheaf theoretical proof of de Rham Theorem).

\bigskip

A large supply of forms in $\Omega_{\lconf}^{k}(M)$ is provided by the following

\begin{proposition}
 Given a locally bounded\footnote{The map is locally bounded if it maps  compact sets to bounded sets.} map $f : M \to N$ in the class $W^{1,n}_{loc}$ between two Riemannian manifolds and a  smooth differential form $\omega \in \Omega_{\mathrm{de Rham}}^{\bullet}(N)$, we have $f^*(\omega) \in  
 \Omega_{\lconf}^{\bullet}(M)$. Furthermore $df^*(\omega) = f^*(d\omega)$. In other words we have a well defined morphism of chain complexes
 $$
  f^* : \Omega_{\mathrm{de Rham}}^{\bullet}(N) \to  \Omega_{\lconf}^{\bullet}(M).
 $$
\end{proposition}

This follows from Proposition \ref{prop.chainrule1}.

\qed

\bigskip

The local  conformal de Rham complex, can be defined without using any smooth structure. Indeed, by Theorem
\ref{th.qcinvariance},  it is enough if the
manifold has a quasiconformal  structure, i.e. an $\r^n$ valued atlas whose transition between local coordinates are quasiconformal maps of domains in $\r^n$.  Sullivan has proved in \cite{sullivan78}  that  any topological $n$-manifold admits a unique quasiconformal structure if  $n\ne 4$, it follows that on any topological manifold of dimension  $n\ne 4$, we may represent cohomology classes by differential forms in
the  local  conformal de Rham complex.

On the other hand, Sullivan and Donaldson have proved that there are topological 4-manifolds which do not carry any quasiconformal structure, and also that  there are  pairs of  quasiconformal   4-manifolds which are homeomorphic but not quasiconformally equivalent. Such exotic  quasiconformal structures even exist on $\bold R^4$, see \cite{DS}.

Vladimir Gol'dshtein \newline
Department of Mathematics, \newline
Ben Gurion University of the Negev, \newline
P.O.Box 653, Beer Sheva, \newline
84105, Israel. \newline
vladimir@bgu.ac.il \newline 
 \vspace{1cm}

Marc Troyanov \newline
Institut de G{\'e}om{\'e}trie, alg{\`e}bre et topologie (IGAT) \newline
B{\^a}timent BCH \\
 \'Ecole Polytechnique F{\'e}derale de Lausanne, \newline
 1015 Lausanne - Switzerland \newline
marc.troyanov@epfl.ch \newline

AMS subjclass: 58A10, 42B, 42B20, 58A14. \newline

\end{document}